\documentstyle[twoside,amstex]{article}
\input amssym.def
\input amssym.tex
\def\bc{\begin{center}}
\def\ec{\end{center}}
\def\no{\noindent}
\def\hang{\hangindent\parindent}
\def\textindent#1{\indent\llap{[#1]\enspace}\ignorespaces}
\def\re{\par\hang\textindent}
\begin{document}
\thispagestyle{empty} \vspace*{3 true cm} \pagestyle{myheadings}
\markboth {\hfill {\sl Huanyin Chen and Marjan Sheibani}\hfill}
{\hfill{\sl UNIQUELY WEAKLY NIL-CLEAN CONDITIONS ON ZERO-DIVISORS}\hfill} \vspace*{-1.5 true cm} \bc{\large\bf UNIQUELY WEAKLY NIL-CLEAN CONDITIONS ON ZERO-DIVISORS}\ec

\vskip6mm
\bc{{\bf Huanyin Chen}\\[2mm]
Department of Mathematics, Hangzhou Normal University\\
Hangzhou 310036, China\\
huanyinchen@@aliyun.com}\ec

\bc{{\bf Marjan
Sheibani}\\[2mm]
Faculty of Mathematics, Statistics and Computer S
cience\\
Semnan University, Semnan, Iran\\
m.sheibani1@@gmail.com}\ec

\vskip10mm
\begin{minipage}{120mm}
\no {\bf Abstract:} An element in a ring $R$ is called uniquely weakly nil-clean if every element in $R$ can be uniquely written as a sum or a difference
of a nilpotent and an idempotent in the sense of very idempotents. The structure of the ring in which every zero-divisor is uniquely weakly nil-clean is completely determined.
We prove that every zero-divisor in a ring $R$ is uniquely weakly nil-clean if and only if $R$ is a D-ring, or $R$ is
abelian, periodic, and $R/J(R)$ is isomorphic to a field $F$, ${\Bbb
Z}_{3}\oplus {\Bbb Z}_{3}$, ${\Bbb Z}_{3}\oplus B$ where $B$
is Boolean, or a Boolean ring. As a specific case, rings in which every zero-divisor $a$ or $-a$ is a nilpotent or an idempotent are also considered. Furthermore, we prove that every zero-divisor in a ring $R$ is uniquely nil-clean if and only if $R$ is a D-ring, or
$R$ is abelian, periodic; and $R/J(R)$ is Boolean.\vskip3mm \no {\bf Key words}: Zero-divisor; Uniquely weakly nil-clean ring;
Uniquely nil-clean ring.\vskip3mm \no {\bf MR(2010) Subject Classification}:
16S34, 16U60, 16U99, 16E50.
\end{minipage}

\vskip15mm \bc{\bf 1. INTRODUCTION}\ec

\vskip4mm \no A ring $R$ is clean provided that every element in $R$ is the sum of a unit and an idempotent. Over the last ten to fifteen years there has been an explosion of interest in this class of rings as well as the many generalizations and variations (cf. [7]). In [12], Diesl introduced a interesting subclass of clean rings.
He called a ring $R$ is nil-clean provided that every element in $R$ is the sum of an nilpotent and an idempotent. If the decomposition is unique, we say that $R$ is uniquely nil-clean. The structure of such rings are very attractive (cf. [3], [8], [12] and [14-15]).
Very recently, Danchev and McGovern investigated those rings in
which every element in $R$ is either a sum or a difference
of a nilpotent and an idempotent. They called such ring is a weakly nil-clean ring. Many properties of commutative weakly nil-clean rings have been obtained in [6] and [10]. The motivation is to extend weakly nil-cleanness to noncommutative rings. We say that $e\in R$ is a very idempotent if $e$ or $-e$ is an idempotent. An element $a\in R$ is weakly nil-clean provided that it can be written as the
sum of a nilpotent and a very idempotent. A weakly nil-clean $a\in R$ is uniquely weakly nil-clean provided that $a=w_1+e_1=w_2+e_2$ where $w_1,w_2$ are nilpotent and $e_1,e_2$ are very idempotents $\Longrightarrow e_1^2=e_2^2$. Notice that in the ring ${\Bbb Z}_4 = {\Bbb Z}/4{\Bbb Z}, -1 = 1+2 =-1+0$,
thus the presentation of $-1$ is not unique in the sense of idempotents. However ${\Bbb Z}_4$ is weakly uniquely nil-clean.

We say that $a\in R$ is a zero-divisor if there exist nonzero
$b,c\in R$ such that $ab=0=ca$. Zero-divisors
occur in many classes of rings. In this article, we are concern on rings in which
every zero-divisor is uniquely weakly nil-clean.
The structure of such rings are completely determined. Furthermore, rings in which every zero-divisor is uniquely nil-clean are also considered.

We call a ring $R$ is uniquely weakly D-nil-clean provided that every zero-divisor in $R$ is
uniquely weakly nil-clean. A ring $R$ is called a
periodic ring if for any $a\in R$ there exist distinct $m,n\in
{\Bbb N}$ such that $a^m=a^n$. A ring $R$ is called a D-ring if every
zero-divisor in $R$ is nilpotent (cf. [1]).
We prove that a ring $R$ is uniquely weakly D-nil-clean if and only if $R$ is a D-ring, or $R$ is
abelian, periodic and $R/J(R)$ is isomorphic to a field $F$, ${\Bbb
Z}_{3}\oplus {\Bbb Z}_{3}$, ${\Bbb Z}_{3}\oplus B$ where $B$
is Boolean, or a Boolean ring. As a specific case, we shall explore rings in which every zero-divisor is a nilpotent or
a very idempotent. A ring $R$ is called uniquely D-nil-clean provided that every zero-divisor in $R$ is uniquely nil-clean.
Moreover, we prove that a
ring $R$ is uniquely D-nil-clean if and only if $R$ is a D-ring, or
$R$ is abelian, periodic; and $R/J(R)$ is Boolean.

Throughout, all rings are associative with an identity. We use
$Id(R), N(R)$ and $J(R)$ to denote the sets of all idempotents,
all nilpotent elements and the Jacobson radical of a ring $R$.
$Z(R)$ and $NZ(R)$ stand for the sets of all zero-divisors and non
zero-divisors of a ring $R$.

\vskip10mm\bc{\bf 2. UNIQUELY WEAKLY NIL-CLEAN RINGS}\ec

\vskip4mm The aim of this is to investigate characterizations of uniquely weakly nil-clean rings which will be repeatedly used in the sequel.
The necessary and sufficient conditions under
which a group ring is uniquely weakly nil-clean are thereby obtained. We begin with

\vskip4mm \hspace{-1.8em} {\bf Lemma 2.1 [2, Theorem 2.28].}\ \
{\it Let $R$ be a ring. Then every element in $R$ is a very
idempotent if and only if $R$ is isomorphic to one of the
following:}
\begin{enumerate}
\item [(a)]{\it ${\Bbb Z}_3$,}
\item [(b)]{\it a Boolean ring, or}
\item [(c)]{\it ${\Bbb Z}_{3}\oplus B$ where $B$ is a Boolean.}
\end{enumerate}

\vskip4mm \hspace{-1.8em} {\bf Theorem 2.2.}\ \ {\it Let $R$ be a
ring. Then $R$ is uniquely weakly nil-clean if and only if} \vspace{-.5mm}
\begin{enumerate}
\item [(1)]{\it $R$ is abelian;}
\vspace{-.5mm}
\item [(2)]{\it $R$ is periodic;}
\vspace{-.5mm}
\item [(3)]{\it $R/J(R)$ is isomorphic to one of the following:}
\begin{enumerate}
\item [(a)]{\it ${\Bbb Z}_3$,}
\item [(b)]{\it a Boolean ring, or}
\item [(c)]{\it ${\Bbb Z}_{3}\oplus B$ where $B$ is a Boolean.}
\end{enumerate}
\end{enumerate}
\vspace{-.5mm} {\it Proof.}\ \ Suppose that $R$ is uniquely weakly nil-clean.
For any idempotent $e\in R$ and any $a\in R$, $e+ea(1-e)\in R$ is
an idempotent. Since $\big(e+ea(1-e)\big)+0=e+ea(1-e)$, by the
uniqueness, $\big(e+ea(1-e)\big)^2=e^2$; hence, $ea(1-e)=0$. This yields that $ea=eae$. Likewise, $ae=eae$, and
so $ea=ae$. Thus, $R$ is abelian. Let $a\in R$. Then there exists a central very idempotent
$e\in R$ such that $w:=a-e\in N(R)$. If $e^2=e$, then
$a-a^2=w-2ew-w^2\in N(R)$. If $e^2=-e$, then $a+a^2=w+2ew+w^2\in
N(R)$. In any case, we can find some $n\in {\Bbb N}$ such that
$a^n=a^{n+1}f(a)$ where $f(t)\in R[t]$. In view of Herstein's
Theorem, $R$ is periodic, and then $N(R)$ forms an ideal of $R$.
Therefore, $J(R)=N(R)$, and so every element in $R/J(R)$ is a very
idempotent. In light of Lemma 2.1, $(3)$ is satisfied.

Conversely, assume that $(1)-(3)$ hold. Let $a\in R$. Then
$\overline{a}$ is a very idempotent, in terms of Lemma 2.1. As
$J(R)$ is nil, every idempotent lifts modulo $J(R)$, and so every
very idempotent lifts modulo $J(R)$. Thus, we can find a very
idempotent $e\in R$ such that $\overline{a}=\overline{e}$. Hence,
$v:=a-e\in J(R)\subseteq N(R)$. If there exists a very idempotent
$f\in R$ such that $w:=a-f\in N(R)$, then
$e^2-f^2=(a-v)^2-(a-w)^2=(-av-va+v^2)+(aw+wa-w^2)$. As $v\in
J(R)$, we see that $-av-va+v^2\in J(R)$. Furthermore,
$aw+wa-w^2\in N(R)$ since $aw=wa$. This implies that
$1-(e^2-f^2)=-(-av-va+v^2)+\big(1-(aw+wa-w^2)\big)\in U(R)$. As
$e^2,f^2\in R$ are idempotents, we have $(e^2-f^2)^3=e^2-f^2$, and
so $(e^2-f^2)\big(1-(e^2-f^2)\big)=0$. Accordingly, $e^2=f^2$, as
asserted.\hfill$\Box$

\vskip4mm \hspace{-1.8em} {\bf Corollary 2.3.}\ \ {\it Let $R$ be a
ring. Then $R$ is uniquely weakly nil-clean if and only if} \vspace{-.5mm}
\begin{enumerate}
\item [(1)]{\it $R$ is abelian;}
\vspace{-.5mm}
\item [(2)]{\it $J(R)$ is nil;}
\vspace{-.5mm}
\item [(3)]{\it $R/J(R)$ is isomorphic to one of the following:}
\begin{enumerate}
\item [(a)]{\it ${\Bbb Z}_3$,}
\item [(b)]{\it a Boolean ring, or}
\item [(c)]{\it ${\Bbb Z}_{3}\oplus B$ where $B$ is a Boolean.}
\end{enumerate}
\end{enumerate}
\vspace{-.5mm} {\it Proof.}\ \ $\Longrightarrow$ In view of Theorem 2.2, $R$ is periodic. Thus, $J(R)$ is nil, as required.

$\Longleftarrow$ By $(3)$, every element in $R/J(R)$ is a very idempotent. By $(2)$, every idempotent lifts modulo $J(R)$. Let $a\in R$. Then $a-a^2\in N(R)$.
As in the proof of Theorem 2.2, $R$ is periodic. This completes the proof, by Theorem 2.2.\hfill$\Box$

\vskip4mm For a local ring $R$, we further derive that $R$ is uniquely weakly nil-clean if and only if
$J(R)$ is nil; $R/J(R)$ is isomorphic to ${\Bbb Z_2}$ or ${\Bbb Z}_3$.

\vskip4mm \hspace{-1.8em} {\bf Corollary 2.4.}\ \ {\it Let $R$ be a
ring. Then $R$ is uniquely weakly nil-clean if and only if}
\begin{enumerate}
\item [(1)]{\it $R$ is periodic;}
\item [(2)]{\it $R$ is uniquely weakly D-nil-clean;}
\item [(3)]{\it $U(R)=\{ x\pm 1~|~x\in N(R)\}$.}
\end{enumerate}
\vspace{-.5mm} {\it Proof.}\ \ Suppose that $R$ is uniquely weakly nil-clean.
In view of Theorem 2.2, $R$ is periodic. $(2)$ is obvious.
Let $x\in U(R)$. Then we have a very idempotent $e\in R$ such that
$w:=x-e\in N(R)$. As $R$ is abelian, we see that $e=x-w$ and
$ew=we$, and so $e=\pm 1$. Therefore $x=w\pm 1$, as desired.

Conversely, assume that $(1)-(3)$ hold. Let $a\in R$. Then we have
distinct $m,n\in {\Bbb N} (m>n)$ such that $a^m=a^n$. If $a$ is a
zero-divisor, then $a$ is uniquely weakly nil-clean. If $a$ is a non
zero-divisor, $a^{m-n}=1$. By $(3)$, we see that $a$ is uniquely weakly nil-clean. This completes the proof.\hfill$\Box$

\vskip4mm Let $P(R)$ be the intersection of all prime ideals of
$R$, i.e., $P(R)$ is the prime radical of $R$. As is well known,
$P(R)$ is the intersection of all minimal prime ideals of $R$.

\vskip4mm \hspace{-1.8em} {\bf Proposition 2.5.}\ \ {\it Let $R$ be a
ring. Then $R$ is uniquely weakly nil-clean if and only if} \vspace{-.5mm}
\begin{enumerate}
\item [(1)]{\it $R$ is abelian;}
\vspace{-.5mm}
\item [(2)]{\it $R/P(R)$ is uniquely weakly nil-clean.}
\end{enumerate}
\vspace{-.5mm} {\it Proof.}\ \ Suppose that $R$ is uniquely weakly nil-clean.
Then $R$ is abelian. In view of Theorem 2.2, $R$ is clean, and so
it is an exchange ring. Thus, $R/P(R)$ is abelian. Obviously,
$J\big(R/P(R)\big)=J(R)/P(R)$ is nil. Further,
$\big(R/P(R)\big)/J\big(R/P(R)\big)\cong R/J(R)$. By Theorem 2.2
again, $R/P(R)$ is uniquely weakly nil-clean.

Conversely, assume that $(1)$ and $(2)$ hold. For any $x\in J(R)$,
we see that $\overline{x}\in J\big(R/P(R)\big)$ is nilpotent.
Since $P(R)$ is nil, we see that $x\in R$ is an nilpotent; hence
that $J(R)$ is nil. As $R/J(R)\cong
\big(R/P(R)\big)/J\big(R/P(R)\big)$, it follows from Theorem 2.2
that $R$ is uniquely weakly nil-clean, as asserted. \hfill$\Box$

\vskip4mm Let $R$ be a ring, and let $G$ be a group. The
augmentation ideal $I(R,G)$ of the group ring $RG$ is the kernel
of the homomorphism from $RG$ to $R$ induced by collapsing $G$ to
$1$. That is, $I(R,G)=ker(\omega)$, where $\omega=\{
\sum\limits_{g\in G}r_gg~|~\sum\limits_{g\in G}r_g=0\}$.

\vskip4mm \hspace{-1.8em} {\bf Lemma 2.6.}\ \ {\it Let $R$ be a
ring, and let $G$ be a group. If $RG$ is uniquely weakly nil-clean, then so
is $R$.} \vskip2mm\hspace{-1.8em} {\it Proof.}\ \ Let $a\in R$.
Then we have a very idempotent $e\in RG$ such that $a-e\in N(RG)$
and that such representation is unique. Hence,
$a-\omega(e)=\omega(a-e)\in N(R)$. Obviously, $\omega(e)\in R$ is
a very idempotent. If $a-f\in N(R)$ for a very idempotent $f\in
R$, then $e=f$, as desired.\hfill$\Box$

\vskip4mm \hspace{-1.8em} {\bf Theorem 2.7.}\ \ {\it Let $R$ be a
ring, and let $G$ be a group. If $I(R,G)$ is nil, then $RG$ is
uniquely weakly nil-clean if and only if so is $R$.} \vskip2mm\hspace{-1.8em}
{\it Proof.}\ \ One direction is obvious by Lemma 2.6. Conversely,
assume that $R$ is uniquely weakly nil-clean. Let $x\in RG$. Then
$x=\omega(x)+\big(x-\omega(x)\big)$. By hypothesis, there exists a
very idempotent $e\in R$ such that $w:=\omega(x)-e\in N(R)$.
Hence, $x=e+\big(w+(x-\omega(x))\big)$. Since $ker(\omega)$ is
nil, we see that $v:=w+(x-\omega(x))\in N(R)$. Assume that $x=f+w$
where $f\in RG$ is an very idempotent and $w\in N(RG)$. Then
$f-\omega(f)\in ker(\omega)$ is nil. As $R$ is uniquely weakly nil-clean, $R$
is abelian. Hence,
$\big(f-\omega(f)\big)\big(1-(f-\omega(f))^2\big)=0$, and so
$f=\omega(f)\in R$. It is easy to verify that
$vw=(x-e)(x-f)=(x-f)(x-e)=wv$, and then $e-f=w-v\in N(R)$. It
follows from $(e-f)\big(1-(e-f)^2\big)=0$ that $e=f$, as
needed.\hfill$\Box$

\vskip4mm \hspace{-1.8em} {\bf Corollary 2.8.}\ \ {\it Let $R$ be
a ring with a prime $p\in J(R)$, and let $G$ be a locally finite
$p$-group. Then $RG$ is uniquely weakly nil-clean if and only if $R$ is uniquely weakly nil-clean.} \vskip2mm\hspace{-1.8em} {\it Proof.}\ \
One direction is obvious. Conversely, assume that $R$ is uniquely weakly nil-clean. Then $J(R)$ is nil by Corollary 2.3. We first suppose $G$ is
finite and prove the claim by induction on $|G|$. As the center of
a nontrivial finite $p$-group contains more than one element, we
may take $x\in G$ be an element in the center with the order $p$.
Let $(x)$ be the subgroup of $G$ generated by $x$. Then
$\overline{G}=G/(x)$ has smaller order. By induction hypothesis,
$ker\big(\overline{\omega}\big)$ is nil, where $\overline{\omega}:
R\overline{G}\to R, \sum\limits_{\overline{g}\in
\overline{G}}r_{\overline{g}}\overline{g}$. Let $\varphi: RG\to
R\overline{G}, \sum\limits_{g}r_gg\to
\sum\limits_{g}r_g\overline{g}$. Then $ker(\varphi)=(1-x)RG$.
Since $x^p=1$, we see that $(1-x)^p\in pRG$; hence, $1-x\in RG$ is
nilpotent. But
$\varphi\big(ker(\omega)\big)=ker\big(\overline{\omega}\big)$ is
nil. For any $z\in ker(\omega)$, we have some $m\in {\Bbb N}$ such
that $z^m\in ker(\varphi)$ is nilpotent. Thus, $z\in RG$ is
nilpotent. We conclude that $ker(\omega)$ is nil, and therefore
$RG$ is uniquely weakly nil-clean, in terms of Theorem 2.7.\hfill$\Box$

\vskip10mm\bc{\bf 3. FACTORIZATION OF ZERO-DIVISORS}\ec \vskip4mm
In this section, we work out the structure of uniquely weakly D-nil-clean rings. To do this, we
need the connections between the class of uniquely weakly nil-clean rings and the class of uniquely weakly D-nil-clean rings.

\vskip4mm \hspace{-1.8em} {\bf Lemma 3.1.}\ \ {\it Every uniquely weakly
D-nil-clean ring is abelian.}\vskip2mm\hspace{-1.8em} {\it Proof.}\
\ Let $e\in R$ be an idempotent, and let $x\in R$. Then
$e+ex(1-e)\in R$ is an idempotent. If $e=1$, then $ex=exe$. If
$1-e=ex(1-e)$, then $ex=exe$. If $e\neq 1$ and $1-e\neq ex(1-e)$,
then $e+ex(1-e)\in R$ is a zero-divisor, as
$$(1-e)\big(e+ex(1-e)\big)=0=\big(e+ex(1-e)\big)\big(1-e-ex(1-e)\big).$$
Since $e+ex(1-e)=e+ex(1-e)+0$, by hypothesis,
$e^2=\big(e+ex(1-e)\big)^2$, and then $ex(1-e)=0$. That is,
$ex=exe$. Likewise, $xe=exe$. Thus, $ex=xe$. This completes the
proof.\hfill$\Box$

\vskip4mm \hspace{-1.8em} {\bf Theorem 3.2.}\ \ {\it Every
uniquely weakly D-nil-clean ring is a D-ring or the product of two uniquely weakly nil-clean rings.}\vskip2mm\hspace{-1.8em} {\it Proof.}\ \ Let $R$
be a uniquely weakly D-nil-clean ring. In view of Lemma 3.1, $R$ is
abelian.

Case I. $R$ is indecomposable. Then every zero-divisor is
nilpotent or invertible. The later is impossible, and so $R$ is a
D-ring.

Case II. $R$ is decomposable. Write $R=A\oplus B$. Let $a\in A$.
Then $(a,0)\in R$ is a zero-divisor. By hypothesis, there exists a very idempotent
$(e,e')\in R$ such that $(a,0)-(e,e')\in N(R)$,
and that $(a,0)-(f,f')\in N(R)$ with a very idempotent $(f,f')\in
R$ implies that $(e,e')^2=(f,f')^2$. Thus, $a-e\in N(R)$. If there
exists a very idempotent $g\in A$ such that $a-g\in N(A)$. Then
$(a,0)-(g,0)\in N(R)$. This implies that $(g,0)^2=(e,e')^2$, and
so $g^2=e^2$. Therefore $A$ is uniquely weakly nil-clean. Similarly, $B$ is
uniquely weakly nil-clean, as asserted.\hfill$\Box$

\vskip4mm \hspace{-1.8em} {\bf Lemma 3.3.}\ \ {\it Let $R$ be a
ring. Then every zero-divisor in $R$ is a very idempotent if and
only if $R$ is isomorphic to one of the following:}
\begin{enumerate}
\item [(1)]{\it a domain,}
\item [(2)]{\it ${\Bbb Z}_{3}\oplus {\Bbb Z}_{3}$,}
\item [(3)]{\it ${\Bbb Z}_{3}\oplus B$ where $B$ is a Boolean, or}
\item [(4)]{\it a Boolean ring.}
\end{enumerate}
\vspace{-.5mm} {\it Proof.}\ \ Suppose that every zero-divisor in
$R$ is a very idempotent. By Lemma 3.1, $R$ is abelian.

Case I. $R$ is indecomposable. Then $Id(R)=\{ 0,1\}$ and
$-Id(R)=\{ 0,-1\}$. Thus, the only zero-divisor is zero. Hence,
$R$ is a domain.

Case II. $R$ is decomposable. Then we have $S,T\neq 0$ such that
$R=S\oplus T$. For any $t\in T$, $(0,t)\in R$ is a zero-divisor.
By hypothesis, $(0,t)$ or $-(0,t)$ is an idempotent; hence that
$t$ or $-t$ is an idempotent in $T$. Therefore every element in
$T$ is a very idempotent. In light of Lemma 2.1, $T$ is isomorphic
to one of the following:
\begin{enumerate}
\item [(a)]{\it ${\Bbb Z}_3$,}
\item [(b)]{\it a Boolean ring, or}
\item [(c)]{\it ${\Bbb Z}_{3}\oplus B$ where $B$ is a Boolean.}
\end{enumerate} Likewise, $S$ is isomorphic to one of the
preceding. Thus, $R$ is isomorphic to one of the following: $R$ is
isomorphic to one of the following:
\begin{enumerate}
\item [(i)]{\it ${\Bbb Z}_3\oplus {\Bbb Z}_3$,}
\item [(ii)]{\it a Boolean ring, or}
\item [(iii)]{\it ${\Bbb Z}_{3}\oplus B$ where $B$ is a Boolean.}
\item [(iv)]{\it ${\Bbb Z}_{3}\oplus {\Bbb Z}_{3}\oplus B$ where $B$ is a Boolean.}
\end{enumerate}
But in Case $(iv)$, $(1,2,0)\in {\Bbb Z}_{3}\oplus {\Bbb
Z}_{3}\oplus B$ is a zero-divisor, while it is not a very
idempotent. Therefore Case $(iv)$ will not appear, as desired.

Conversely, if $R$ is a domain, then every zero-divisor is zero.
If $R={\Bbb Z}_{3}\oplus {\Bbb Z}_{3}$, then $NZ(R)=\{
(1,1),(1,2),(2,1)\}, Id(R)=\{ (0,0),(0,1),(1,0),(1,1)\}$ and
$-Id(R)=\{ (0,0),$ $(0,2),$ $(2,0),$ $(2,2)\}$. Therefore
$R=NZ(R)\bigcup Id(R)\bigcup -Id(R)$. If $R={\Bbb Z}_{3}\bigoplus
B$ where $B$ is a Boolean, then $Id(R)=\{ (0,b),(1,b)~|~b\in B\}$
and $-Id(R)=\{ (0,b), (2,b)~|~b\in B\}$. Therefore $R=Id(R)\bigcup
-Id(R)$. If $R$ is a Boolean ring, then every element in $R$ is an
idempotent. In any case, every element in $R$ is a very
idempotent, and we are done.\hfill$\Box$

\vskip4mm We come now to the main result of this section.

\vskip4mm \hspace{-1.8em} {\bf Theorem 3.4.}\ \ {\it Let $R$ be a
ring. Then $R$ is uniquely weakly D-nil-clean if and only if $R$ is a D-ring, or }
\vspace{-.5mm}
\begin{enumerate}
\item [(1)]{\it $R$ is abelian;}
\vspace{-.5mm}
\item [(2)]{\it $R$ is periodic;}
\vspace{-.5mm}
\item [(3)]{\it $R/J(R)$ is isomorphic to one of the following:}
\begin{enumerate}
\item [(a)]{\it a field $F$,}
\item [(b)]{\it ${\Bbb Z}_{3}\oplus {\Bbb Z}_{3}$,}
\item [(c)]{\it ${\Bbb Z}_{3}\oplus B$ where $B$ is Boolean, or}
\item [(d)]{\it a Boolean ring.}
\end{enumerate}
\end{enumerate}
\vspace{-.5mm} {\it Proof.}\ \ $\Longrightarrow $ Suppose that $R$ is not a D-ring. In view of Theorem 3.2, $R$ is the product of two uniquely
weakly nil-clean rings $R_1$ and $R_2$. By Theorem 2.2, $R_1$ and $R_2$ are abelian periodic rings, and then so is $R$.
In view of [4, Theorem], $N(R)$ is an
ideal of $R$. As $R$ is periodic, $J(R)$ is nil; hence,
$J(R)=N(R)$. As every idempotent lifts modulo $N(R)$, we see that
$R/J(R)$ is abelian. Let $\overline{a}\in R/N(R)$ be a
zero-divisor. If $a\in R$ is not a zero-divisor, then $a\in U(R)$,
and so $\overline{a}\in U\big(R/N(R)\big)$, a contradiction. Thus,
$a\in R$ is a zero-divisor. By hypothesis, $a$ is the sum of a
very idempotent and a nilpotent. Hence, $\overline{a}$ is a very
idempotent. That is, every zero-divisor in $R/J(R)$ is a very
idempotent. In light of Lemma 3.3, $R/J(R)$ is isomorphic to one
of the following:
\begin{enumerate}
\item [(i)]{\it a domain $F$,}
\item [(ii)]{\it ${\Bbb Z}_{3}\oplus {\Bbb Z}_{3}$,}
\item [(iii)]{\it ${\Bbb Z}_{3}\oplus B$ where $B$ is Boolean, or}
\item [(iv)]{\it a Boolean ring.}
\end{enumerate} If $R=F$ is a domain, then for any $a\in R$, $a=0$ or $a^m=1$ for some $m\in {\Bbb N}$.
This shows that $R$ is a field, as required.

$\Longleftarrow $ In view of [4,
Theorem ], $N(R)$ forms an ideal of $R$. Let $a\in R$ be a zero
divisor. Then $\overline{a}\in R/J(R)$ is a zero-divisor;
otherwise, $\overline{a}\in R/J(R)$ is invertible, and so $a\in R$
is invertible, a contradiction. According to Lemma 3.5,
$\overline{a}$ is a very idempotent in $R/J(R)$. As $R$ is
periodic, $J(R)$ is nilpotent, and so every idempotent modulo
$J(R)$. This implies that $v:=a-e\in N(R)$ for some very
idempotent $e\in R$. Let $f\in R$ be a very idempotent such that
$w:=a-f\in N(R)$. Then
$e^2-f^2=(a-v)^2-(a-w)^2=-av-va+v^2+aw+wa-w^2\in N(R)$. As $e,f\in
R$ are very clean, we se that $e^2,f^2\in R$ are idempotents. It
is easy to verify that $(e^2-f^2)\big(1-(e^2-f^2)^2\big)=0$, and
so $e^2=f^2$. Therefore we complete the proof.\hfill$\Box$\hfill$\Box$

\vskip4mm We now consider a specific case and explore the structure
of rings in which every zero-divisor is a very idempotent or a
nilpotent element.

\vskip4mm \hspace{-1.8em} {\bf Lemma 3.5.}\ \ {\it Every ring in
which every element is a very idempotent or a nilpotent element is
abelian.} \vskip2mm\hspace{-1.8em} {\it Proof.}\ \ Let $e\in R$ be
an idempotent, and let $x\in R$. Then $1-ex(1-e)\in U(R)$. If
$\big(1-ex(1-e)\big)^2=1-ex(1-e)$, then $ex(1-e)=0$, and so
$ex=exe$. If $\big(1-ex(1-e)\big)^2=-\big(1-ex(1-e)\big),$ then
$ex(1-e)=2$. and so $ex(1-e)=2e(1-e)=0$. Hence, $ex=exe$. If
$1-ex(1-e)\in N(R)$, this will be a contradiction. Thus, $ex=exe$.
Likewise, $xe=exe$. Therefore $ex=xe$, hence the
result.\hfill$\Box$

\vskip4mm \hspace{-1.8em} {\bf Lemma 3.6.}\ \ {\it Let $R$ be a
ring. Then the following are equivalent:}
\begin{enumerate}
\item [(1)]{\it $R=N(R)\cup Id(R) \cup -Id(R)$.}
\item [(2)]{\it $R=J(R)\cup Id(R)\cup -Id(R)$}
\item [(3)]{\it $R$ is isomorphic to one of the following:}
\end{enumerate}\begin{enumerate}
\item [(a)]{\it ${\Bbb Z}_3$, ${\Bbb Z}_4$,}
\item [(b)]{\it a Boolean ring, or}
\item [(c)]{\it ${\Bbb Z}_{3}\oplus B$ where $B$ is a Boolean.}
\end{enumerate}
\vspace{-.5mm} {\it Proof.}\ \ $(1)\Leftrightarrow $ This is proved by [10, Proposition 1.21]

$(1)\Leftrightarrow$ This is obvious by [10, Proposition 1.19] and Lemma 2.1.\hfill$\Box$

\vskip4mm \hspace{-1.8em} {\bf Theorem 3.7.}\ \ {\it Let $R$ be a
ring. Then $R$ is an abelian ring in which every zero-divisor in
$R$ is a very idempotent or a nilpotent element if and only if $R$
is isomorphic to one of the following:}
\begin{enumerate}
\item [(a)]{\it a D-ring,}
\item [(b)]{\it a Boolean ring,}
\item [(c)]{\it ${\Bbb Z}_3\oplus {\Bbb Z}_3$,}
\item [(d)]{\it ${\Bbb Z}_{3}\oplus B$ where $B$ is a Boolean.}
\end{enumerate}
\vspace{-.5mm} {\it Proof.}\ \ Suppose that $R$ is an abelian ring
in which every zero-divisor in $R$ is a very idempotent or a
nilpotent element.

Case I. $R$ is indecomposable. Then every very idempotent is $0,1$
or $-1$. Hence, every zero-divisor in $R$ is nilpotent. Hence, $R$
is a D-ring.

Case II. $R$ is decomposable. Write $R=S\oplus T$. For any $t\in
T$, $(0,t)$ is a very idempotent or a nilpotent element. We infer
that every element in $T$ is a very idempotent or a nilpotent
element. Similarly, every element in $S$ is a very idempotent or a
nilpotent element. By virtue of Lemma 3.6, $S$ and $T$ are both
isomorphic to one of the following:
\begin{enumerate}
\item [(a)]{\it ${\Bbb Z}_3$, ${\Bbb Z}_4$,}
\item [(b)]{\it a Boolean ring, or}
\item [(c)]{\it ${\Bbb Z}_{3}\oplus B$ where $B$ is a Boolean.}
\end{enumerate}
But one easily checks that $Z(R)\neq Id(R)\bigcup -Id(R)\bigcup
N(R)$ for any of those types
\begin{enumerate}
\item [(1)]{\it ${\Bbb Z}_3\oplus {\Bbb Z}_4$,}
\item [(2)]{\it ${\Bbb Z}_3\oplus {\Bbb Z}_3\oplus B$ where $B$ is a Boolean ring,}
\item [(3)]{\it ${\Bbb Z}_{4}\oplus {\Bbb Z}_4$,}
\item [(4)]{\it ${\Bbb Z}_4\oplus B$ where $B$ is a Boolean ring, and }
\item [(5)]{\it ${\Bbb Z}_{3}\oplus {\Bbb Z}_4\oplus B$ where $B$ is a Boolean ring.}
\end{enumerate} Therefore $R$ is isomorphic
to one of $(a)-(d)$.

Conversely, $R$ is abelian, as every D-ring is connected. One
easily checks that any of these four types of rings satisfy the
desired condition.\hfill$\Box$

\vskip4mm By Theorem 3.7 and Theorem 3.4, we see that every abelian ring in which every zero-divisor in
$R$ is a very idempotent or a nilpotent element is uniquely weakly D-nil-clean.

\vskip4mm \hspace{-1.8em} {\bf Corollary 3.8.}\ \ {\it Let $R$ be
a ring. Then the following are equivalent:}
\vspace{-.5mm}
\begin{enumerate}
\item [(1)]{\it $R$ is an abelian ring in which every zero-divisor in $R$ is an idempotent or a nilpotent
element.} \vspace{-.5mm}
\item [(2)]{\it $R$ is a D-ring or a Boolean ring.}
\end{enumerate}
\vspace{-.5mm} {\it Proof.}\ \ $(1)\Rightarrow (2)$ In view of
Theorem 3.7, $R$ is isomorphic to one of the following:
\begin{enumerate}
\item [(a)]{\it a D-ring,}
\item [(b)]{\it a Boolean ring,}
\item [(c)]{\it ${\Bbb Z}_3\oplus {\Bbb Z}_3$,}
\item [(d)]{\it ${\Bbb Z}_{3}\oplus B$ where $B$ is a Boolean.}
\end{enumerate} But in the case ${\Bbb Z}_3\oplus {\Bbb Z}_3$,
$(0,2)\not\in Id(R)\bigcup N(R)$. In the case ${\Bbb
Z}_{3}\oplus B$, $(2,0)\not\in Id(R)\bigcup N(R)$. Therefore
proving $(2)$.

$(2)\Rightarrow (1)$ This is obvious.\hfill$\Box$

\vskip4mm The abelian condition
in Corollary 3.8 is necessary. For instance, every zero-divisor
in $T_2\big({\Bbb Z}_2\big)$ is an idempotent or a nilpotent
element. But $T_2\big({\Bbb Z}_2\big)$ is neither a Boolean ring
nor a D-ring.

\vskip10mm\bc{\bf 4. UNIQUELY D-NIL-CLEAN RINGS}\ec

\vskip4mm The aim of this is to give the connection of uniquely D-nil-clean rings and
uniquely nil-clean rings, and then characterize the structure of
uniquely D-nil-clean rings.

\vskip4mm \hspace{-1.8em} {\bf Lemma 4.1.}\ \ {\it Every
uniquely D-nil-clean ring is abelian.} \vskip2mm\hspace{-1.8em}
{\it Proof.}\ \ This is similar to that in Lemma 3.1.\hfill$\Box$

\vskip4mm \hspace{-1.8em} {\bf Proposition 4.2.}\ \ {\it A ring
$R$ is uniquely D-nil-clean if and only if for any zero-divisor
$a\in R$ there exists a central idempotent $e\in R$ such that
$a-e\in N(R)$.}\vskip2mm\hspace{-1.8em} {\it Proof.}\ \ One
direction is obvious from Lemma 4.1. Conversely, letting $e\in R$
be an idempotent, we have a central idempotent $f\in R$ such that
$w:=e-f\in N(R)$. Thus, $(e-f)^3=e-f$, and so
$(e-f)\big(1-(e-f)^2\big)=0$. This implies that $e=f$, and then
$R$ is abelian. Let $a\in R$ be a zero-divisor. Then there exists
a central $e\in R$ such that $a-e\in N(R)$. If there exists an
idempotent $f\in R$ such that $a-f\in N(R)$, then
$e-f=(a-f)-(a-e)\in N(R)$. It follows from $(e-f)^3=e-f$ that
$e=f$, which completes the proof.\hfill$\Box$

\vskip4mm \hspace{-1.8em} {\bf Lemma 4.3.}\ \ {\it Let $R$ be a
ring. Then $R$ is a uniquely D-nil-clean ring if and only if $R$
is a D-ring or $R$ is uniquely nil-clean.}
\vskip2mm\hspace{-1.8em} {\it Proof.}\ \ $\Longrightarrow$ In view of Lemma 4.1, $R$ is abelian.

Case I. $R$ is indecomposable. Let $a\in R$ be a zero-divisor.
Then $a\in R$ is nilpotent or $a\in U(R)$. This shows that every
zero-divisor is nilpotent, i.e., $R$ is a D-ring.

Case II. $R$ is decomposable. Write $R=A\oplus B$. For any $x\in
A$, $(x,0)\in R$ is a zero-divisor. Hence, we can find a unique
idempotent $(e,f)\in R$ such that $(x,0)-(e,f)\in N(R)$. Thus,
$x-e\in N(R)$ for an idempotent $e\in R$. If there exists an
idempotent $g\in R$ such that $x-g\in N(R)$. Then $(x,0)-(g,f)\in
N(R)$. By the uniqueness, we get $g=e$, and therefore $A$ is
uniquely nil-clean. Similarly, $B$ is uniquely nil-clean, and then
$R$ is uniquely nil-clean.

$\Longleftarrow$ If $R$ is a D-ring, then $R$ is a D-uniquely nil
clean ring. So we assume that $R$ is uniquely nil-clean, and
therefore $R$ is a uniquely D-nil-clean ring.\hfill$\Box$

\vskip4mm \hspace{-1.8em} {\bf Theorem 4.4.}\ \ {\it Let $R$ be a ring. Then $R$ is uniquely D-nil-clean if and only if $R$ is a D-ring, or}
\vspace{-.5mm}
\begin{enumerate}
\item [(1)]{\it $R$ is abelian;}
\vspace{-.5mm}
\item [(2)]{\it $R$ is periodic;}
\vspace{-.5mm}
\item [(3)]{\it $R/J(R)$ is Boolean.}
\end{enumerate}
\vspace{-.5mm} {\it Proof.}\ \ $\Longrightarrow$ In view of Lemma 4.1, $R$ is abelian. Suppose that $R$ is not a D-ring. Then $R$ is a uniquely nil-clean
ring, in terms of Lemma 4.3. In view of [12, Theorem 5.9], $R/J(R)$ is Boolean and $J(R)$ is nil. Let $a\in R$. Then $a-a^2\in J(R)$, and so $(a-a^2)^m=0$ for some $m\in {\Bbb N}$. Similarly to Theorem 2.2, $R$ is periodic, in terms of Herstein¡¯s Theorem.

In view of [4, Theorem], $N(R)$ forms an ideal of $R$. Hence,
$J(R)=N(R)$. Let $\overline{a}\in R/J(R)$ is a zero-divisor. Then
$a\in R$ is a divisor; otherwise, $a\in U(R)$ as $R$ is periodic,
a contradiction. Hence, $a$ is the sum of an idempotent and a
nilpotent element. This shows that $\overline{a}$ is an
idempotent. Therefore, every zero-divisor in $R/J(R)$ is an
idempotent.

Set $S=R/J(R)$. Suppose that $S$ has a nonzero zero-divisor. Then
we have some $x,y\in R$ such that $xy=0, x,y\neq 0$. Hence,
$(yx)^2=0$. If $yx\neq 0$, then $yx\in R$ is a zero-divisor. So
$yx\in R$ is an idempotent. Thus, $yx=(yx)^2=0$. This implies that
$x\in R$ is a zero-divisor, and so $x=x^2$. It follows that
$1-x\in R$ is a zero-divisor; hence that $1-x=(1-x)^2$. Therefore
$x^2=x$.

Let $a\in R$. Then $\big(xa(1-x)\big)^2=0$. Hence, $xa(1-x)=0$;
otherwise, $xa(1-x)\in R$ is an idempotent, and so $xa(1-x)=0$, a
contradiction. Thus, $xa(1-x)=0$, hence, $xa=xax$. Likewise,
$ax=xax$. Thus, $xa=ax$. If $xa=0$, then $a\in R$ is a
zero-divisor, and so it is an idempotent. If $xa\neq 0$, it
follows from $xa(1-x)=0$ that $xa\in R$ is a zero-divisor, and so
$xa=(xa)^2$. Hence, $xa(1-a)=0$. This implies that $1-a\in R$ is a
zero-divisor, and then $1-a=(1-a)^2$. Thus, $a=a^2$. Therefore
$a\in R$ is an idempotent. Consequently, $R/J(R)$ is Boolean or
$R/J(R)$ is a domain. If $R/J(R)$ is a domain, the periodic
property implies that $R$ is a field. Thus, $R$ is local. But $J(R)$ is nil, and so every zero-divisor is nilpotent. We infer that $R$ is a D-ring, an absurd.
This shows that $R/J(R)$ is Boolean, as desired.

$\Longleftarrow$ If $R$ is a D-ring, then $R$ is uniquely D-nil-clean. We now assume that $(1)-(3)$ hold.
Let $a\in R$ be a
zero-divisor. As $R/J(R)$ is Boolean, $a-a^2\in J(R)\subseteq N(R)$. Thus, we can find an idempotent
$e\in R$ such that $a-e\in N(R)$. Since $R$ is
abelian, we see that such idempotent $e$ is unique. Therefore $R$
is uniquely D-nil-clean.\hfill$\Box$

\vskip4mm By Theorem 4.4 and Theorem 3.4, every uniquely D-nil-clean ring is a uniquely weakly D-nil-clean ring.

\vskip4mm \hspace{-1.8em} {\bf Lemma 4.5.}\ \ {\it Let $R$ be a
ring. Then $R$ is uniquely nil-clean if and only if}\vspace{-.5mm}
\begin{enumerate}
\item [(1)]{\it $R$ is abelian.}
\vspace{-.5mm}
\item [(2)]{\it $R/J(R)$ is Boolean and $J(R)$ is nil.}
\vspace{-.5mm}
\end{enumerate}\vspace{-.5mm}{\it Proof.}\ \ $\Longrightarrow $ This is obvious by [12, Lemma 5.5 and Theorem 5.9].

$\Longleftarrow$ For any $a\in R$, $a-a^2\in
J(R)$, and so we have an idempotent $e\in R$ such that $a-e\in
J(R)$, as $J(R)$ is nil. Write $a=e+v$. Then $v\in J(R)\subseteq
N(R)$. If there exists an idempotent $f\in R$ and a $w\in N(R)$
such that $a=f+w$, then $e-f=(a-v)-(a-w)=w-v$. Clearly,
$wv=(a-f)(a-e)=(a-e)(a-f)=vw$, and so $e-f\in N(R)$. Since
$(e-f)^3=e-f$, we see that $e-f=0$, and then $e=f$. Therefore $R$
is uniquely nil clean. \hfill$\Box$

\vskip4mm \hspace{-1.8em} {\bf Theorem 4.6.}\ \ {\it Let $R$ be a
ring. Then $R$ is uniquely nil-clean if and only if}
\vspace{-.5mm}
\begin{enumerate}
\item [(1)]{\it $2\in R$ is nilpotent;}
\vspace{-.5mm}
\item [(2)]{\it $R$ is uniquely weakly nil-clean.}
\end{enumerate}
\vspace{-.5mm} {\it Proof.}\ \ Suppose that $R$ is uniquely nil
clean. In view of Lemma 4.5, $\overline{2}^2=\overline{2}$ in
$R/J(R)$, and so $2\in J(R)$ is nilpotent. By Lemma 4.5 and
Theorem 2.2, we observe that every uniquely nil-clean ring is uniquely weakly nil-clean.

Conversely, assume that $(1)$ and $(2)$ hold. As $2\in {\Bbb Z}_3$
is not nilpotent. In view of Theorem 2.2, $R$ is
abelian, $J(R)$ is nil, and that $R/J(R)$ is Boolean. The result
follows by Lemma 4.5.\hfill$\Box$

\vskip4mm \hspace{-1.8em} {\bf Corollary 4.7.}\ \ {\it Let $R$ be
a ring. Then $R$ is uniquely nil-clean if and only if}
\vspace{-.5mm}
\begin{enumerate}
\item [(1)]{\it $R$ is abelian;}
\vspace{-.5mm}
\item [(2)]{\it $R/P(R)$ is uniquely nil-clean.}
\end{enumerate}
\vspace{-.5mm} {\it Proof.}\ \ One direction is obvious, by
Theorem 4.6 and Proposition 2.5.

Conversely, assume that $(1)$ and $(2)$ hold. By virtue of Theorem
4.6, $\overline{2}\in R/P(R)$ is nilpotent. We infer that $2\in R$
is nilpotent. Furthermore, $R/P(R)$ is uniquely weakly nil-clean. According
to Proposition 2.5, $R$ is uniquely weakly nil-clean. By using Theorem 4.6 again,
$R$ is uniquely nil-clean.\hfill$\Box$

\vskip4mm \hspace{-1.8em} {\bf Corollary 4.8.}\ \ {\it Let $R$ be
a ring, and $G$ be a group. Then $RG$ is uniquely nil-clean if and
only if $R$ is uniquely nil-clean and $I(R,G)$ is nil.}
\vskip2mm\hspace{-1.8em} {\it Proof.}\ \ Suppose $RG$ is uniquely
nil-clean. Then $RG$ is uniquely weakly nil-clean and $2\in N(RG)$, by
Theorem 4.6. Hence, $R$ is uniquely weakly nil-clean and $2\in N(R)$. By
using Theorem 4.6 again, $R$ is uniquely nil-clean. Thanks to Lemma 4.5,
$RG/J(RG)$ is Boolean. For any $g\in G$, we see that
$(1-g)-(1-g)^2\in J(RG)$; hence, $1-g\in J(RG)$. This implies that
$ker(\omega)\subseteq J(RG)$ is nil, as desired.

Conversely, assume that $R$ is uniquely nil-clean and
$ker(\omega)$ is nil. In light of Theorem 4.6 and Theorem 2.7,
$2\in N(R)$ and $RG$ is uniquely weakly nil-clean. By using Theorem 4.6
again, $RG$ is uniquely nil-clean.\hfill$\Box$

\vskip4mm Let $G$ be a
$3$-group. Then ${\Bbb Z}_3G$ is is not uniquely nil-clean by Corollary 4.8, while
it is uniquely weakly nil-clean.

\vskip4mm \hspace{-1.8em} {\bf Corollary 4.9.}\ \ {\it Let $R$ be a ring in which $2$ is nonnilpotent. Then $R$ is uniquely $D$-nil-clean if and only if $R$ is uniquely weakly D-nil-clean.}\vskip2mm\hspace{-1.8em} {\it Proof.}\ \ $\Longrightarrow$ This is obvious.

$\Longleftarrow $ In light of Theorem 3.2, $R$ is a D-ring, or the product of two uniquely weakly nil-clean rings. As $2$ is nonnilpotent in $R$,
every $R$ is uniquely weakly nil-clean is uniquely nil-clean, by Theorem 4.6. Therefore $R$ is uniquely $D$-nil-clean, in terms of Lemma 4.3.
 \hfill$\Box$

\vskip15mm \bc{\bf REFERENCES}\ec \vskip4mm {\small \re{1} H.
Abu-Khuzam and A. Yaqub, Structure of rings with certain
conditions on zero-divisors, {\it Int. J. Math. Math. Sci.}, {\bf
15}(2006), 1--6.

\re{2} M.S. Ahn, Weakly Clean Rings and Almost Clean rings, Ph.D.
Thesis, The University of Iowa, Iowa, 2003.

\re{3} D. Andrica and G. C\u{a}lug\u{a}reanu, A nil-clean $2\times
2$ matrix over the integers which is not clean, {\it J. Algebra
Appl.}, {\bf 13}(2014), 1450009 [9 pages] DOI:
10.1142/S0219498814500091.

\re{4} A. Badawi, On abelian $\pi$--regular rings, {\it{Comm.
Algebra,}} {\bf{25}}(1997), 1009--1021.

\re{5} S. Breaz; G. Galugareanu; P.
Danchev and T. Micu, Nil-clean matrix rings, {\it Linear Algebra
Appl.}, {\bf 439}(2013), 3115--3119.

\re{6} S. Breaz; P. Danchev and Y.
Zhou, Weakly nil-clean rings, preprint, 2015.

\re{7} H.
Chen, {\it Rings Related Stable Range Conditions}, Series in
Algebra 11, World Scientific, Hackensack, NJ, 2011.

\re{8} H. Chen, Strongly nil-clean matrices over
$R[x]/\big(x^2-1\big)$, {\it Bull. Korean Math. Soc.}, {\bf
49}(2012), 589--599.

\re{9} H. Chen, On strongly nil-clean matrices, {\it Comm.
Algebra}, {\bf 41}(2013), 1074--1086.

\re{10} P.V. Danchev and W.W. McGovern,
Commutative weakly nil-clean unital rings, {\it J. Algebra}, {\bf 425}(2015), 410--422.

\re{11} A.J. Diesl, Classes of Strongly
Clean Rings, Ph.D. Thesis, University of California, Berkeley,
2006.

\re{12} A.J. Diesl, Nil clean rings,
 {\it J. Algebra}, {\bf 383}(2013), 197--211.

\re{13} V.A. Hiremath and S. Hegde, Using ideals to provide a
unified approach to uniquely clean rings, {\it J. Aust. Math.
Soc.}, {\bf 96}(2014), 258--274.

\re{14} M.T. Kosan; T.K. Lee and Y. Zhou, When is every matrix
over a division ring a sum of an idempotent and a nilpotent? {\it
Linear Algebra Appl.}, {\bf 450}(2014), 7--12.

\re{15} A. Stancu, On some constuctions of nil-clean, clean, and
exchange rings, arXiv: 1404.2662v1 [math.RA], 10 Apr 2014.
\end{document}